\documentclass[12pt]{article}

\usepackage{amsmath, epsfig, cite}
\usepackage{amssymb}
\usepackage{amsfonts}
\usepackage{latexsym}
\usepackage{float}

\newtheorem{thm}{Theorem}[section]

\newtheorem{cor}[thm]{Corollary}

\newtheorem{lem}[thm]{Lemma}

\numberwithin{equation}{section}

\newcommand{\qed}{{\hfill$\square$}\medskip}
\newcommand*\pFqskip{8mu}
\catcode`,\active
\newcommand*\pFq{\begingroup
        \catcode`\,\active
        \def ,{\mskip\pFqskip\relax}%
        \dopFq
}
\catcode`\,12
\def\dopFq#1#2#3#4#5{%
        {}_{#1}F_{#2}\biggl[\genfrac..{0pt}{}{#3}{#4};#5\biggr]%
        \endgroup
}

\begin{document}


\begin{center}
{\Large\bf A $p$-adic supercongruence for truncated\\[10pt]
hypergeometric series ${}_7F_6$}
\end{center}

\vskip 2mm \centerline{Ji-Cai Liu}
\begin{center}
{\footnotesize College of Mathematics and Information Science, Wenzhou University, Wenzhou 325035, PR China\\
{\tt jcliu2016@gmail.com, jc2051@163.com} }
\end{center}


\vskip 0.7cm \noindent{\bf Abstract.}
Using an identity due to Gessel and Stanton and some properties of the $p$-adic Gamma function, we establish a $p$-adic supercongruence for truncated hypergeometric series ${}_7F_6$. From it we deduce some related supercongruences, which extend certain recent results and confirm a supercongruence conjecture.

\vskip 2mm
\noindent{\it MR Subject Classifications}: Primary 11A07; Secondary 33C20

\section{Introduction}
Throughout this paper, let $p$ denote an odd prime. Based on a result of Ahlgren and Ono \cite{AO}, Kilbourn proved\cite{Kilbourn} that
\begin{align}
\sum_{k=0}^{\frac{p-1}{2}}\left(\frac{\left(\frac{1}{2}\right)_k}{k!}\right)^4
\equiv a_p \pmod{p^3},\label{aa1}
\end{align}
where $(x)_k=x(x+1)\cdots (x+k-1)$ and $a_p$ is the $p$-th coefficient of a weight $4$ modular form
\begin{align*}
\eta(2z)^4\eta(4z)^4:=q\prod_{n\ge 1}(1-q^{2n})^4(1-q^{4n})^4,\quad q=e^{2\pi iz}.
\end{align*}
For more such supercongruences, one refers to Rodriguez-Villegas \cite{RV}.

Motivated by Ramanujan-type formulas for $1/\pi$, Zudilin \cite{Zudilin} obtained several supercongruences by using the WZ method. For example, he proved that
\begin{align}
\sum_{k=0}^{\frac{p-1}{2}}(-1)^k(4k+1)\left(\frac{\left(\frac{1}{2}\right)_k}{k!}\right)^3
\equiv (-1)^{\frac{p-1}{2}}p\pmod{p^3},\label{aa2}
\end{align}
which was conjectured by van Hamme \cite{vanHamme} and first confirmed by Mortenson \cite{Mortenson}.

McCarthy and Osburn \cite{MO} proved another conjecture of van Hamme \cite{vanHamme}:
\begin{align}
\sum_{k=0}^{\frac{p-1}{2}}(-1)^k(4k+1)\left(\frac{\left(\frac{1}{2}\right)_k}{k!}\right)^5
\equiv \begin{cases}
-\frac{p}{\Gamma_p\left(\frac{3}{4}\right)^4} \pmod{p^3}\quad &\text{if $p\equiv 1 \pmod{4}$,}\\
0\pmod{p^3}\quad &\text{if $p\equiv 3 \pmod{4}$,}
\end{cases}\label{aa3}
\end{align}
where $\Gamma_p(\cdot)$ stands for the $p$-adic Gamma function recalled in next section.

Hypergeometric series are very important in many research fields, including algebraic varieties, differential equations and modular forms. It is known that some of the truncated hypergeometric series are related to the number of rational points on certain algebraic varieties over finite fields and further to coefficients of modular forms.
For complex numbers $a_i, b_j$ and $z$ with none of the $b_j$ being negative integers or zero, we define the truncated hypergeometric series as
\begin{align*}
\pFq{r}{s}{a_1,a_2,\cdots,a_r}{b_1,b_2,\cdots,b_s}{z}_n=\sum_{k=0}^{n}\frac{(a_1)_k (a_2)_k\cdots (a_r)_k}{(b_1)_k (b_2)_k\cdots (b_s)_k}\cdot \frac{z^k}{k!}.
\end{align*}
All of the left-hand sides of \eqref{aa1}-\eqref{aa3} can be rewritten as the truncated hypergeometric series.

Using the Dougall's formula and some properties of the $p$-adic Gamma function,
Long and Ramakrishna \cite{LR} proved that for any prime $p\ge 5$,
\begin{align*}
\pFq{7}{6}{\frac{1}{3},\frac{1}{3},\frac{1}{3},\frac{1}{3},\frac{1}{3},\frac{1}{3},\frac{7}{6}}{1,1,1,1,1,\frac{1}{6}}{1}_{p-1}\equiv \begin{cases}
 -p\Gamma_p\left(\frac{1}{3}\right)^9 \pmod{p^6}\quad &\text{if $p\equiv 1\pmod{6}$,}\\[10pt]
-\frac{10}{27}p^4\Gamma_p\left(\frac{1}{3}\right)^9 \pmod{p^6}\quad &\text{if $p\equiv 5\pmod{6}$,}
\end{cases}
\end{align*}
which extends a result of van Hamme \cite{vanHamme}.

Some other interesting supercongruences for the truncated hypergeometric series were obtained by several authors, see, for example, \cite{BS,He,Long,LR,Mortenson2005}.

Let $\langle x\rangle_p$ denote the least non-negative integer $r$ with $x\equiv r \pmod{p}$ and
$\lfloor x \rfloor$ denote the greatest integer less than or equal to a real number $x$. The aim of this paper is to establish the following supercongruence.

\begin{thm}\label{tt1}
Let $p\ge 5$ be a prime. For any $p$-adic integer $\alpha$ with $0\le \langle \alpha \rangle_p\le \lfloor p/4 \rfloor$, we have
\begin{align}
&\pFq{7}{6}{\frac{1}{2},\frac{1}{2},\frac{1}{2},\frac{1}{4},\frac{7}{6},\frac{1}{2}+\alpha,\frac{1}{4}-\alpha}
{1,1,1,\frac{1}{6},1+2\alpha,\frac{1}{2}-2\alpha}{1}_{\frac{p-1}{2}}\notag\\[10pt]
&\equiv\begin{cases}
(-1)^{\frac{p+3}{4}}p\Gamma_p\left(\frac{1}{2}\right)\Gamma_p\left(\frac{1}{4}\right)^2
\Gamma_p\left(1+\alpha\right)\Gamma_p\left(\frac{3}{4}-\alpha\right)\\[5pt] \times\Gamma_p\left(\frac{1}{2}+\alpha\right)^3\Gamma_p\left(\frac{1}{4}-\alpha\right)^3
\pmod{p^3}\quad &\text{if $p\equiv 1\pmod{4}$},\\[10pt]
0\pmod{p^3}\quad &\text{if $p\equiv 3\pmod{4}$}.
\end{cases}\label{aa4}
\end{align}
\end{thm}

Some related supercongruences for truncated hypergeometric series can be deduced from
\eqref{aa4}.
\begin{cor}\label{tt2}
For any prime $p\ge 5$, we have
\begin{align}
&\pFq{5}{4}{\frac{1}{2},\frac{1}{2},\frac{1}{2},\frac{1}{4},\frac{7}{6}}{1,1,1,\frac{1}{6}}{\frac{1}{4}}_{\frac{p-1}{2}}\notag\\
&\equiv \begin{cases} (-1)^{\frac{p+3}{4}}p\Gamma_p\left(\frac{1}{2}\right)\Gamma_p\left(\frac{1}{4}\right)^2 \pmod{p^3}\quad &\text{if $p\equiv 1\pmod{4}$,}\\[10pt]
0\pmod{p^3}\quad &\text{if $p\equiv 3\pmod{4}$.}
\end{cases}\label{aa5}
\end{align}
\end{cor}

\begin{cor}\label{tt3}
For any prime $p\ge 5$, we have
\begin{align}
&\pFq{6}{5}{\frac{1}{2},\frac{1}{2},\frac{1}{2},\frac{1}{4},\frac{1}{4},\frac{7}{6}}{1,1,1,1,\frac{1}{6}}{1}_{\frac{p-1}{2}}\notag\\
&\equiv \begin{cases}
-p\Gamma_p\left(\frac{1}{4}\right)^4 \pmod{p^3}\quad &\text{if $p\equiv 1\pmod{4}$,}\\[10pt]
0\pmod{p^3}\quad &\text{if $p\equiv 3\pmod{4}$.}
\end{cases}\label{aa6}
\end{align}
\end{cor}

Supercongruences \eqref{aa5} and \eqref{aa6} extend some recent results of He \cite[Theorem 1.1]{He}, which stated the corresponding supercongruences hold mod $p^2$. It should be mentioned that the factor $(-1)^{\frac{p+3}{4}}$ was missing in the first case of \cite[(1.1)]{He}. Supercongruence \eqref{aa5} confirms the second conjectural supercongruence in \cite[Conjecture 1.2]{He}.

In the next section, we first recall some properties of the Morita's $p$-adic Gamma function and a terminating hypergeometric series identity due to Gessel and Stanton. The proof of Theorem \ref{tt1} and Corollary \ref{tt2} and \ref{tt3} will be given in Section 3.

\section{Some lemmas}
We first recall some basic properties of the Morita's $p$-adic Gamma function. For more details, one refers to \cite[\S 11.6]{Cohen} and \cite{Morita}.
Let $p$ be an odd prime and $\mathbb{Z}_p$ denote the set of all $p$-adic integers. For
$x\in \mathbb{Z}_p$, the Morita's $p$-adic Gamma function \cite[Definition 11.6.5]{Cohen} is defined as
\begin{align*}
\Gamma_p(x)=\lim_{m\to x}(-1)^m\prod_{\substack{0\le k < m\\
(k,p)=1}}k,
\end{align*}
where the limit is for $m$ tending to $x$ $p$-adically in $\mathbb{Z}_{\ge 0}$.
\begin{lem}
Suppose $p$ is an odd prime and $x\in \mathbb{Z}_p$. Then
\begin{align}
&\Gamma_p(1)=-1,\label{cc1}\\
&\Gamma_p(x)\Gamma_p(1-x)=(-1)^{a_p(x)},\label{cc2-old}\\
&\frac{\Gamma_p(x+1)}{\Gamma_p(x)}=
\begin{cases}
-x\quad&\text{if $v_p(x)=0$,}\\
-1\quad &\text{if $v_p(x)>0$, }
\end{cases}\label{cc3}
\end{align}
where $a_p(x)\in \{1,2,\cdots,p\}$ with $x\equiv a_p(x)\pmod{p}$ and $v_p(\cdot)$ denotes the $p$-order.
\end{lem}

In fact, \eqref{cc2-old} can be extended as follows.
\begin{lem} (Long and Ramakrishna \cite[(2.3)]{LR})
Let $p\ge 5$ be a prime, $y\in \mathbb{C}_p$ and $x\in \mathbb{Q}$ with $v_p(x)\ge 0$. Then
\begin{align}
\Gamma_p(x+y)\Gamma_p(1-x-y)=(-1)^{a_p(x)}. \label{cc2}
\end{align}
\end{lem}

\begin{lem} (Long and Ramakrishna \cite[Lemma 17, (4)]{LR})
Let $p$ be an odd prime. If $a\in \mathbb{Z}_p, n\in \mathbb{N}$
such that none of $a,a+1,\cdots,a+n-1$ in $p\mathbb{Z}_p$, then
\begin{align}
(a)_n=(-1)^n\frac{\Gamma_p(a+n)}{\Gamma_p(a)}.\label{cc4}
\end{align}
\end{lem}

The following lemma is a special case of the theorem due to Long and Ramakrishna \cite[Theorem 14]{LR}.
\begin{lem}
Suppose $p\ge 5$ is a prime. If $a\in \mathbb{Z}_p,m\in \mathbb{C}_p$
satisfying $v_p(m)\ge 0$, then
\begin{align}
\Gamma_p(a+mp)\equiv \Gamma_p(a)\sum_{k=0}^2 \frac{G_k(a)}{k!}(mp)^k \pmod{p^3},\label{cc5}
\end{align}
where $G_k(a)=\Gamma_p^{(k)}(a)/\Gamma_p(a)\in \mathbb{Z}_p$ and $\Gamma_p^{(k)}(x)$ is the k-th
derivative of $\Gamma_p(x)$.
\end{lem}

\begin{lem} (Gessel and Stanton \cite[(1.8)]{GS})
If $n$ is a non-negative integer, then
\begin{align}
&\pFq{7}{6}{a,b,a-b+\frac{1}{2},1+\frac{2a}{3},1-2d,2a+2d+n,-n}{2a-2b+1,2b,\frac{2a}{3},a+d+\frac{1}{2}
,1-d-\frac{n}{2},1+a+\frac{n}{2}}{1}\notag\\
&=\begin{cases}
\displaystyle\frac{\left(\frac{1}{2}\right)_r \left(b+d\right)_r \left(d-b+a+\frac{1}{2}\right)_r \left(a+1\right)_r}{\left(b+\frac{1}{2}\right)_r \left(a+d+\frac{1}{2}\right)_r \left(d\right)_r \left(a-b+1\right)_r}\quad &\text{if $n=2r$},\\[15pt]
0\quad &\text{if $n=2r+1$}.
\end{cases}\label{cc6}
\end{align}
\end{lem}

\section{Proof of \eqref{aa4}--\eqref{aa6}}
{\noindent \it Proof of \eqref{aa4}.}
Let $\omega$ be any primitive 3th root of unity.
Letting $n=\frac{p-1}{2},a=\frac{1}{4},b=\frac{1}{2}+\alpha,d=\frac{1+\omega^2p}{4}$ in \eqref{cc6} and noting that $1+\omega+\omega^2=0$, we obtain
\begin{align}
&\pFq{7}{6}{\frac{1-p}{2},\frac{1-\omega p}{2},\frac{1-\omega^2p}{2},\frac{1}{2}+\alpha,\frac{1}{4}-\alpha,\frac{1}{4},\frac{7}{6}}{1+\frac{p}{4},1+\frac{\omega p}{4},1+\frac{\omega^2p}{4},1+2\alpha,\frac{1}{2}-2\alpha,\frac{1}{6}}{1}\notag\\[10pt]
&=\begin{cases}
\displaystyle\frac{\left(\frac{1}{2}\right)_{\frac{p-1}{4}} \left(\frac{5}{4}\right)_{\frac{p-1}{4}} \left(\frac{4\alpha+3+\omega^2p}{4}\right)_{\frac{p-1}{4}} \left(\frac{2-4\alpha+\omega^2p}{4}\right)_{\frac{p-1}{4}}}{\left(\frac{1+\omega^2p}{4}\right)_{\frac{p-1}{4}}\left(\frac{4+\omega^2p}{4}\right)_{\frac{p-1}{4}}(1+\alpha)_{\frac{p-1}{4}}
\left(\frac{3}{4}-\alpha\right)_{\frac{p-1}{4}}}\quad &\text{if $p\equiv 1\pmod{4}$},\\
0\quad &\text{if $p\equiv 3\pmod{4}$}.
\end{cases}\label{ff1}
\end{align}

Since for $0\le k\le \frac{p-1}{2}$ and $0\le \langle\alpha \rangle_p\le \lfloor p/4\rfloor$,
\begin{align*}
\left(1+2\alpha\right)_k\left(\frac{1}{2}-2\alpha\right)_k\equiv
\left(1+2\langle\alpha \rangle_p\right)_k\left(\frac{1}{2}-2\langle\alpha \rangle_p\right)_k
\not\equiv 0\pmod{p},
\end{align*}
and $\left(\frac{7}{6}\right)_k/\left(\frac{1}{6}\right)_k=6k+1$,
we conclude that none of the denominators on the left-hand side of \eqref{ff1}
contain a multiple of $p$.
By the fact that
\begin{align*}
(u+vp)(u+vp\omega)(u+vp\omega^2)=u^3+v^3p^3,
\end{align*}
we have
\begin{align*}
(u+vp)_k(u+vp\omega)_k(u+vp\omega^2)_k\equiv (u)_k^3\pmod{p^3},
\end{align*}
and so
\begin{align}
&\pFq{7}{6}{\frac{1-p}{2},\frac{1-\omega p}{2},\frac{1-\omega^2p}{2},\frac{1}{2}+\alpha,\frac{1}{4}-\alpha,\frac{1}{4},\frac{7}{6}}{1+\frac{p}{4},1+\frac{\omega p}{4},1+\frac{\omega^2p}{4},1+2\alpha,\frac{1}{2}-2\alpha,\frac{1}{6}}{1}\notag\\[5pt]
&\equiv \pFq{7}{6}{\frac{1}{2},\frac{1}{2},\frac{1}{2},\frac{1}{4},\frac{7}{6},\frac{1}{2}+\alpha,\frac{1}{4}-\alpha}
{1,1,1,\frac{1}{6},1+2\alpha,\frac{1}{2}-2\alpha}{1}_{\frac{p-1}{2}}\pmod{p^3}.\label{ff2}
\end{align}

Combining \eqref{ff1} and \eqref{ff2}, we conclude the proof of \eqref{aa4} for $p\equiv 3\pmod{4}$.
In order to prove the case $p\equiv 1 \pmod{4}$, it suffices to show that the following holds mod $p^3$:
\begin{align}
&\frac{\left(\frac{1}{2}\right)_{\frac{p-1}{4}} \left(\frac{5}{4}\right)_{\frac{p-1}{4}} \left(\frac{4\alpha+3+\omega^2p}{4}\right)_{\frac{p-1}{4}} \left(\frac{2-4\alpha+\omega^2p}{4}\right)_{\frac{p-1}{4}}}{\left(\frac{1+\omega^2p}{4}\right)_{\frac{p-1}{4}}\left(\frac{4+\omega^2p}{4}\right)_{\frac{p-1}{4}}(1+\alpha)_{\frac{p-1}{4}}
\left(\frac{3}{4}-\alpha\right)_{\frac{p-1}{4}}}\notag\\
&\equiv (-1)^{\frac{p+3}{4}}p\Gamma_p\left(\frac{1}{2}\right)\Gamma_p\left(\frac{1}{4}\right)^2
\Gamma_p\left(1+\alpha\right)\Gamma_p\left(\frac{3}{4}-\alpha\right)\Gamma_p\left(\frac{1}{2}+\alpha\right)^3\Gamma_p\left(\frac{1}{4}-\alpha\right)^3.
\label{ff3}
\end{align}

By \eqref{cc4}, we have
\begin{align}
 \left(\frac{5}{4}\right)_{\frac{p-1}{4}}=p\left(\frac{1}{4}\right)_{\frac{p-1}{4}}
 =(-1)^{\frac{p-1}{4}}p\frac{\Gamma_p\left(\frac{p}{4}\right)}{\Gamma_p\left(\frac{1}{4}\right)},\label{ff4}
\end{align}
and
\begin{align}
&\frac{\left(\frac{1}{2}\right)_{\frac{p-1}{4}} }{ \left(\frac{1+\omega^2p}{4}\right)_{\frac{p-1}{4}}\left(\frac{4+\omega^2p}{4}\right)_{\frac{p-1}{4}}}\notag\\
&=(-1)^{\frac{3(p-1)}{4}}\frac{\Gamma_p\left(\frac{1+p}{4}\right)\Gamma_p\left(\frac{1+\omega^2p}{4}\right)\Gamma_p\left(1+\frac{\omega^2p}{4}\right)}
{\Gamma_p\left(\frac{1}{2}\right)\Gamma_p\left(\frac{-\omega p}{4}\right)\Gamma_p\left(\frac{3-\omega p}{4}\right)}\notag\\
&=(-1)^{\frac{3(p-1)}{4}}\frac{\Gamma_p\left(\frac{1+p}{4}\right)\Gamma_p\left(\frac{1+\omega p}{4}\right)\Gamma_p\left(\frac{1+\omega^2p}{4}\right)\Gamma_p\left(1+\frac{ p}{4}\right)\Gamma_p\left(1+\frac{\omega p}{4}\right)\Gamma_p\left(1+\frac{\omega^2p}{4}\right)}
{\Gamma_p\left(\frac{1}{2}\right)\Gamma_p\left(1+\frac{ p}{4}\right)\Gamma_p\left(\frac{-\omega p}{4}\right)\Gamma_p\left(1+\frac{\omega p}{4}\right)\Gamma_p\left(\frac{3-\omega p}{4}\right)\Gamma_p\left(\frac{1+\omega p}{4}\right)\label{ff5}
},
\end{align}
where we have utilized the fact that $1+\omega+\omega^2=0$ in the first step.
Applying \eqref{cc5} and the symmetry with respect to the 3th roots of unity, we get
\begin{align}
&\Gamma_p\left(\frac{1+p}{4}\right)\Gamma_p\left(\frac{1+\omega p}{4}\right)\Gamma_p\left(\frac{1+\omega^2p}{4}\right)\Gamma_p\left(1+\frac{ p}{4}\right)\Gamma_p\left(1+\frac{\omega p}{4}\right)\Gamma_p\left(1+\frac{\omega^2p}{4}\right)\notag\\
&\equiv \left(\Gamma_p(1)\Gamma_p\left(\frac{1}{4}\right)\right)^3\pmod{p^3} \notag\\
&= -\Gamma_p\left(\frac{1}{4}\right)^3.\quad \text{(by \eqref{cc1})}\label{ff6}
\end{align}
Furthermore, by \eqref{cc2} we get
\begin{align}
&\Gamma_p\left(\frac{-\omega p}{4}\right)\Gamma_p\left(1+\frac{\omega p}{4}\right)\Gamma_p\left(\frac{3-\omega p}{4}\right)\Gamma_p\left(\frac{1+\omega p}{4}\right)\notag\\
&=(-1)^{a_p\left(1\right)+a_p\left(\frac{1}{4}\right)}\notag\\
&=(-1)^{1+\frac{3p+1}{4}}\notag\\
&=(-1)^{\frac{3(p-1)}{4}}.\label{ff7}
\end{align}
Finally, combining \eqref{ff4}--\eqref{ff7} we obtain
\begin{align}
\frac{\left(\frac{1}{2}\right)_{\frac{p-1}{4}} \left(\frac{5}{4}\right)_{\frac{p-1}{4}} }{ \left(\frac{1+\omega^2p}{4}\right)_{\frac{p-1}{4}}\left(\frac{4+\omega^2p}{4}\right)_{\frac{p-1}{4}}}
&\equiv (-1)^{\frac{p+3}{4}}p\frac{\Gamma_p\left(\frac{1}{4}\right)^2\Gamma_p\left(\frac{p}{4}\right)}
{\Gamma_p\left(\frac{1}{2}\right)\Gamma_p\left(1+\frac{p}{4}\right)}\pmod{p^3}.\label{ff8}
\end{align}
From \eqref{cc2} and \eqref{cc3}, we have
\begin{align}
\frac{\Gamma_p\left(\frac{p}{4}\right)}{\Gamma_p\left(1+\frac{p}{4}\right)}=-1,\label{ff9}
\end{align}
and
\begin{align}
\Gamma_p\left(\frac{1}{2}\right)^2=(-1)^{\frac{p+1}{2}}=-1.\label{ff10}
\end{align}
Substituting \eqref{ff9} and \eqref{ff10} into \eqref{ff8} gives
\begin{align}
\frac{\left(\frac{1}{2}\right)_{\frac{p-1}{4}} \left(\frac{5}{4}\right)_{\frac{p-1}{4}} }{ \left(\frac{1+\omega^2p}{4}\right)_{\frac{p-1}{4}}\left(\frac{4+\omega^2p}{4}\right)_{\frac{p-1}{4}}}
\equiv (-1)^{\frac{p+3}{4}}p\Gamma_p\left(\frac{1}{2}\right)\Gamma_p\left(\frac{1}{4}\right)^2 \pmod{p^3}.\label{ff11}
\end{align}

It is not hard to verify that for $0\le \langle \alpha \rangle_p\le \lfloor p/4\rfloor$,
none of
\begin{align*}
\left(\frac{4\alpha+3+\omega^2p}{4}\right)_{\frac{p-1}{4}},\quad
\left(\frac{2-4\alpha+\omega^2p}{4}\right)_{\frac{p-1}{4}},\quad
(1+\alpha)_{\frac{p-1}{4}},\quad
\left(\frac{3}{4}-\alpha\right)_{\frac{p-1}{4}}
\end{align*}
contain a multiple of $p$. It follows from \eqref{cc4} that
\begin{align}
&\frac{\left(\frac{4\alpha+3+\omega^2p}{4}\right)_{\frac{p-1}{4}} \left(\frac{2-4\alpha+\omega^2p}{4}\right)_{\frac{p-1}{4}}}{(1+\alpha)_{\frac{p-1}{4}}
\left(\frac{3}{4}-\alpha\right)_{\frac{p-1}{4}}}\notag\\
&=\frac{\Gamma_p\left(\frac{4\alpha+2-\omega p}{4}\right)\Gamma_p\left(\frac{1-4\alpha-\omega p}{4}\right)\Gamma_p\left(1+\alpha\right)\Gamma_p\left(\frac{3}{4}-\alpha\right)}
{\Gamma_p\left(\frac{4\alpha+3+\omega^2 p}{4}\right)\Gamma_p\left(\frac{2-4\alpha+\omega^2 p}{4}\right)\Gamma_p\left(\frac{4\alpha +3+p}{4}\right)\Gamma_p\left(\frac{2-4\alpha+p}{4}\right)}\notag\\
&=\frac{\Gamma_p\left(\frac{4\alpha+2-\omega p}{4}\right)\Gamma_p\left(\frac{2-4\alpha+\omega p}{4}\right)\Gamma_p\left(\frac{1-4\alpha-\omega p}{4}\right)\Gamma_p\left(\frac{4\alpha +3+\omega p}{4}\right)\Gamma_p\left(1+\alpha\right)\Gamma_p\left(\frac{3}{4}-\alpha\right)}
{\Gamma_p\left(\frac{4\alpha+3+p}{4}\right)\Gamma_p\left(\frac{4\alpha+3+\omega p}{4}\right)
\Gamma_p\left(\frac{4\alpha +3+\omega^2 p}{4}\right)\Gamma_p\left(\frac{2-4\alpha+p}{4}\right)
\Gamma_p\left(\frac{2-4\alpha+\omega p}{4}\right)\Gamma_p\left(\frac{2-4\alpha +\omega^2 p}{4}\right)},\label{ff12}
\end{align}
where we have used the fact that $1+\omega+\omega^2=0$ in the first step.

By \eqref{cc2}, we have
\begin{align}
&\Gamma_p\left(\frac{4\alpha+2-\omega p}{4}\right)\Gamma_p\left(\frac{2-4\alpha+\omega p}{4}\right)\Gamma_p\left(\frac{1-4\alpha-\omega p}{4}\right)\Gamma_p\left(\frac{4\alpha +3+\omega p}{4}\right)\notag\\
&=(-1)^{a_p\left(\frac{1-2\alpha}{2}\right)+a_p\left(\frac{4\alpha+3}{4}\right)}.\label{ff13}
\end{align}
Using \eqref{cc5} and the symmetry with respect to the 3th roots of unity, we get
\begin{align}
&\Gamma_p\left(\frac{4\alpha+3+p}{4}\right)\Gamma_p\left(\frac{4\alpha+3+\omega p}{4}\right)
\Gamma_p\left(\frac{4\alpha +3+\omega^2 p}{4}\right)\notag\\
&\times \Gamma_p\left(\frac{2-4\alpha+p}{4}\right)
\Gamma_p\left(\frac{2-4\alpha+\omega p}{4}\right)\Gamma_p\left(\frac{2-4\alpha +\omega^2 p}{4}\right)\notag\\
&\equiv \Gamma_p\left(\frac{4\alpha+3}{4}\right)^3\Gamma_p\left(\frac{1-2\alpha}{2}\right)^3\pmod{p^3}.\label{ff14}
\end{align}
Substituting \eqref{ff13} and \eqref{ff14} into \eqref{ff12} gives
\begin{align}
&\frac{\left(\frac{4\alpha+3+\omega^2p}{4}\right)_{\frac{p-1}{4}} \left(\frac{2-4\alpha+\omega^2p}{4}\right)_{\frac{p-1}{4}}}{(1+\alpha)_{\frac{p-1}{4}}
\left(\frac{3}{4}-\alpha\right)_{\frac{p-1}{4}}}\notag\\
&\equiv (-1)^{a_p\left(\frac{1-2\alpha}{2}\right)+a_p\left(\frac{4\alpha+3}{4}\right)}
\frac{\Gamma_p\left(1+\alpha\right)\Gamma_p\left(\frac{3}{4}-\alpha\right)}
{\Gamma_p\left(\frac{4\alpha+3}{4}\right)^3\Gamma_p\left(\frac{1-2\alpha}{2}\right)^3}\pmod{p^3}.
\label{ff15}
\end{align}

By \eqref{cc2}, we have
\begin{align}
&\Gamma_p\left(\frac{4\alpha+3}{4}\right)^3\Gamma_p\left(\frac{1-4\alpha}{4}\right)^3=
(-1)^{3a_p\left(\frac{4\alpha+3}{4}\right)},\label{ff16}\\
&\Gamma_p\left(\frac{1-2\alpha}{2}\right)^3\Gamma_p\left(\frac{1+2\alpha}{2}\right)^3=
(-1)^{3a_p\left(\frac{1-2\alpha}{2}\right)}.\label{ff17}
\end{align}
Combining \eqref{ff15}--\eqref{ff17}, we obtain
\begin{align}
&\frac{\left(\frac{4\alpha+3+\omega^2p}{4}\right)_{\frac{p-1}{4}} \left(\frac{2-4\alpha+\omega^2p}{4}\right)_{\frac{p-1}{4}}}{(1+\alpha)_{\frac{p-1}{4}}
\left(\frac{3}{4}-\alpha\right)_{\frac{p-1}{4}}}\notag\\
&\equiv \Gamma_p\left(1+\alpha\right)\Gamma_p\left(\frac{3}{4}-\alpha\right)
\Gamma_p\left(\frac{1}{2}+\alpha\right)^3\Gamma_p\left(\frac{1}{4}-\alpha\right)^3\pmod{p^3}.
\label{ff18}
\end{align}
Then the proof of \eqref{ff3} follows from \eqref{ff11} and \eqref{ff18}.
\qed

{\noindent \it Proof of \eqref{aa5}.}
Letting $\alpha\to \infty$ in \eqref{ff1} and noting that
\begin{align*}
\lim_{\alpha\to \infty}\frac{\left(\frac{1}{2}+\alpha\right)_k\left(\frac{1}{4}-\alpha\right)_k}
{\left(1+2\alpha\right)_k\left(\frac{1}{2}-2\alpha\right)_k}=\left(\frac{1}{4}\right)^k,
\end{align*}
we obtain
\begin{align}
&\pFq{5}{4}{\frac{1-p}{2},\frac{1-\omega p}{2},\frac{1-\omega^2p}{2},\frac{1}{4},\frac{7}{6}}{1+\frac{p}{4},1+\frac{\omega p}{4},1+\frac{\omega^2p}{4},\frac{1}{6}}{\frac{1}{4}}\notag\\[10pt]
&=\begin{cases}
\displaystyle\frac{\left(\frac{1}{2}\right)_{\frac{p-1}{4}} \left(\frac{5}{4}\right)_{\frac{p-1}{4}} }{\left(\frac{1+\omega^2p}{4}\right)_{\frac{p-1}{4}}\left(\frac{4+\omega^2p}{4}\right)_{\frac{p-1}{4}}}\quad &\text{if $p\equiv 1\pmod{4}$},\\
0\quad &\text{if $p\equiv 3\pmod{4}$}.
\end{cases}\label{ff19}
\end{align}
The proof of \eqref{aa5} follows from \eqref{ff11}, \eqref{ff19} and the fact that
\begin{align*}
(u+vp)_k(u+vp\omega)_k(u+vp\omega^2)_k\equiv (u)_k^3\pmod{p^3}.
\end{align*}
\qed

{\noindent \it Proof of \eqref{aa6}.}
Letting $\alpha=0$ in \eqref{aa4} reduces to
\begin{align}
&\pFq{6}{5}{\frac{1}{2},\frac{1}{2},\frac{1}{2},\frac{1}{4},\frac{1}{4},\frac{7}{6}}{1,1,1,1,\frac{1}{6}}{1}_{\frac{p-1}{2}}\notag\\
&\equiv \begin{cases}
(-1)^{\frac{p+3}{4}}p\Gamma_p\left(\frac{1}{2}\right)^4\Gamma_p\left(\frac{1}{4}\right)^5
\Gamma_p\left(\frac{3}{4}\right)\Gamma_p\left(1\right) \pmod{p^3}\quad &\text{if $p\equiv 1\pmod{4}$,}\\[10pt]
0\pmod{p^3}\quad &\text{if $p\equiv 3\pmod{4}$.}\label{ff20}
\end{cases}
\end{align}
For $p\equiv 1 \pmod{4}$, by \eqref{cc1} and \eqref{cc2} we have
\begin{align*}
&\Gamma_p\left(1\right)=-1.\\
&\Gamma_p\left(\frac{1}{2}\right)^4=\left((-1)^{\frac{p+1}{2}}\right)^2=1,\\
&\Gamma_p\left(\frac{1}{4}\right)\Gamma_p\left(\frac{3}{4}\right)=(-1)^{\frac{p+3}{4}}.
\end{align*}
Substituting the above equations into \eqref{ff20}, we complete the proof of \eqref{aa6}.
\qed

\vskip 5mm \noindent{\bf Acknowledgments.} The author would like to thank the referees for their careful reviews and useful comments. The author is especially grateful to Professor Dennis Stanton for providing valuable references related to the identity
\eqref{cc6}, which help him to simplify a previous version of this paper.

\end{document}